\pdfoutput=1
\documentclass[11pt,a4]{article}
\usepackage{fullpage}
\usepackage{amsthm, amsmath, amssymb}
\usepackage{color}

\theoremstyle{definition}
\newtheorem{definition}{Definition}[section]
\theoremstyle{theorem}
\newtheorem{theorem}[definition]{Theorem}
\newtheorem{lemma}[definition]{Lemma}
\newtheorem{proposition}[definition]{Proposition}

\DeclareMathOperator{\PP}{PP}

\newcommand{\id}{\mathrm{id}}

\begin{document}

\title{
On refined enumerations of plane partitions of a given shape with bounded entries
}
\author{
Takuya Inoue
\footnote{
Graduate School of Mathematical Sciences, the University of Tokyo,
3-8-1 Komaba, Meguro-ku, 153-8914 Tokyo, Japan.
Email: inoue@ms.u-tokyo.ac.jp
}
}
\date{}
\maketitle

\begin{abstract}
In this paper, we consider plane partitions $\PP(\lambda; m)$ of a given shape $\lambda$, with entries at most $m$.
We prove that the distributions of two statistics on $\PP(\lambda; m)$ coincide:
one is the number of rows containing $0$ and the other is the number of rows containing $m$.
We also provide a bijective proof.
\end{abstract}

\newcommand{\sijto}{\mathrel{\circ \hspace{-1.8mm} = \hspace{-2.2mm} > \hspace{-2.2mm} = \hspace{-1.8mm} \circ}}

\section{Introduction} \label{sec::intro}
The notion of sijections (signed bijections) and the Garsia-Milne involution principle have been widely used in many papers on combinatorics. (See, for example, \cite{Doyle,FK1}.)
In \cite{Inoue} we established the notion of compatibility of sijections, which acts as a bridge between refined enumerations and bijective proofs and provides a way to assess the ``naturalness'' of sijections.
In this paper, we present further applications of compatibilities by combining them with non-intersecting paths and the LGV lemma.
The following is our main theorem:
\begin{theorem} \label{thm::main}
Let $\lambda$ be a partition. Then, the generating function of $\PP(\lambda; m)$ with respect to the number of rows containing $0$, equals that of $\PP(\lambda; m)$ with respect to the number of rows containing $m$. Namely, we have
\begin{equation}
	\sum_{P \in \PP(\lambda; m)} x^{\text{\# of rows containing $0$}} = \sum_{P \in \PP(\lambda; m)} x^{\text{\# of rows containing $m$}}. \label{eq::main}
\end{equation}
\end{theorem}
Here, $\PP(\lambda; m)$ denotes the set of plane partitions of shape $\lambda$ with entries between $0$ and $m$, inclusive.
In fact, the RHS of (\ref{eq::main}) is a special case of the refined enumerations considered by Krattenthaler in \cite{Kra}, where the formula is presented in a more complicated form than ours. For more details, see Proposition \ref{prop::main}.

After proving the main theorem via an explicit determinant formula in Section \ref{sec::main}, we give a bijective proof for it in Section \ref{sec::bij}.
In Section \ref{sec::DLGV}, we formalize the method used in the bijective proof and illustrate the reason why Schur functions are symmetric, as an application of this method.

\section{Notations \& Preliminaries} \label{sec::prel}
Let $\Gamma$ be a directed acyclic graph and $\mathbf{a}=(a_1,a_2,\ldots,a_n)$ and $\mathbf{b}=(b_1,b_2,\ldots,b_n)$ sequences of vertices of $\Gamma$.
Then, we denote paths from $\mathbf{a}$ to $\mathbf{b}$ by $P_\Gamma(\mathbf{a},\mathbf{b})$ and non-intersecting paths from $\mathbf{a}$ to $\mathbf{b}$ by $P_\Gamma^{\text{NI}}(\mathbf{a},\mathbf{b})$.

We briefly review the definitions and properties of signed sets, sijections and compatibilities. See also \cite{FK1,Inoue}.
First, a signed set is a pair of disjoint finite sets and a sijection from a signed set $S=(S^+,S^-)$ to a signed set $T=(T^+,T^-)$ is a bijection between $S^{+} \sqcup T^{-}$ and $S^{-} \sqcup T^{+}$.
We denote a sijection from $S$ to $T$ by $S \sijto T$.
A \textit{statistic} of a signed set $S=(S^+,S^-)$ is a function on $S^+ \sqcup S^-$.
Disjoint unions and Cartesian products of signed sets, sijections, and statistics are defined in a straightforward manner; for details see \cite[Section 2 and 3]{Inoue}.
Let $S$ and $T$ be signed sets and let $\eta_S$ and $\eta_T$ be a statistic of $S$ and  $T$, respectively.
Then, a sijection $\phi \colon S \sijto T$ is \textit{compatible with $\eta = \eta_S \sqcup \eta_T$} if and only if
\[
	\forall s \in S^+ \sqcup S^- \sqcup T^+ \sqcup T^-,\,\eta(\phi(s))=\eta(s),
\]
and we denote it by $(S,\eta_S) \sijto (T,\eta_T)$.
The most important property of the compatibility is that it is preserved under compositions of sijections:
\begin{lemma}[{\cite[Lemma 3.2]{Inoue}}]\label{lem:comp_comp}
Let $\phi \colon S \sijto T$, $\psi \colon T \sijto U$ be sijections and $\eta$ a statistic of $S \sqcup T \sqcup U$.
If $\phi$ and $\psi$ are compatible with $\eta$, then $\psi \circ \phi$ is compatible with $\eta$.
\end{lemma}
For additional properties and further details, see \cite[Section 3]{Inoue}.

\section{Computational Proof} \label{sec::main}
In this section, we give a computational proof of the main theorem.
\begin{theorem} \label{thm::main}
Let $\lambda$ be a partition. Then, the generating function of $\PP(\lambda; m)$ with respect to the number of rows containing $0$, equals that of $\PP(\lambda; m)$ with respect to the number of rows containing $m$. Namely, we have
\[
	\sum_{P \in \PP(\lambda; m)} x^{\text{\# of rows containing $0$}} = \sum_{P \in \PP(\lambda; m)} x^{\text{\# of rows containing $m$}}.
\]
\end{theorem}


Let $\Gamma = (\mathbb{Z}^2, E)$, where $E := \{ (i,j) \to (i,j-1) \} \cup \{ (i,j) \to (i+1,j) \}$.
Let $a_i = (-i,-i)$ and $b_j=(\lambda_j-j,-m-j)$ for $1 \leq i,j \leq \ell(\lambda)$, where $\ell(\lambda)$ is the number of rows of $\lambda$. Then, the plane partitions $\PP(\lambda; m)$ are in bijection to non-intersecting paths $\mathcal{P}_{\Gamma}^{\text{NI}}(\mathbf{a},\mathbf{b})$,
where $\mathbf{a}=(a_1,a_2,\ldots,a_{\ell(\lambda)})$, $\mathbf{b}=(b_1,b_2,\ldots,b_{\ell(\lambda)})$ and each path records entries in corresponding row.
Under this bijection, the number of rows containing $0$ corresponds to the number of edges in $E_1 := \{ b_j+(0,1) \to b_j\}$ that constitute the corresponding paths,
and the number of rows containing $m$ corresponds to the number of edges in $E_2 := \{ a_i \to a_i+(1,0) \}$ that constitute the corresponding paths.
Therefore, we obtain the following proposition by the LGV lemma.

\begin{proposition} \label{prop::main}
It holds that
\begin{align*}
	\sum_{P \in \PP(\lambda; m)} x^{\text{\# of rows containing $0$}} &=  \det \left( \dbinom{\lambda_j+m-1}{m+j-i}x + \dbinom{\lambda_j+m-1}{m+j-i-1} \right)_{ij},\\
	\sum_{P \in \PP(\lambda; m)} x^{\text{\# of rows containing $m$}} &= \det \left( \dbinom{\lambda_j+m-1}{m+j-i}x + \dbinom{\lambda_j+m-1}{m+j-i-1} \right)_{ij}.
\end{align*}
\end{proposition}
This proposition implies Theorem \ref{thm::main}.

\section{Bijective Proof} \label{sec::bij}
The final part of the proof of Theorem \ref{thm::main} can be carried out in a more combinatorial way, avoiding explicit determinant formula.
As is well known, an involution on $\mathcal{P}_{\Gamma}(\mathbf{a},\mathbf{b}) \setminus \mathcal{P}_{\Gamma}^{\text{NI}}(\mathbf{a},\mathbf{b})$ provides a bijective proof of the LGV lemma.
This involution is given by ordering the intersections and swapping the paths after the first intersection.
By combining this with the identity on $\mathcal{P}_{\Gamma}^{\text{NI}}(\mathbf{a},\mathbf{b})$, we obtain a sijection $\Phi_{\text{LGV}} \colon \mathcal{P}_{\Gamma}^{\text{NI}}(\mathbf{a},\mathbf{b}) \sijto \mathcal{P}_{\Gamma}(\mathbf{a},\mathbf{b})$.
Let $\eta_t$ be a statistic on $\mathcal{P}_{\Gamma}(\mathbf{a},\mathbf{b})$ such that $\eta_t(P)$ represents the number of edges in $E_t$ that constitute $P$. 
Then, the sijection $\Phi_{\text{LGV}}$ is compatible with $\eta_1$ and $\eta_2$ since it does not affect the multiset of edges that constitute the paths.
Therefore, if we have a compatible sijection $\Psi \colon (\mathcal{P}_{\Gamma}(\mathbf{a},\mathbf{b}),\eta_1) \sijto (\mathcal{P}_{\Gamma}(\mathbf{a},\mathbf{b}),\eta_2)$,
then we obtain Theorem \ref{thm::main} from the following diagram (and the bijection between $\PP(\lambda; m)$ and $\mathcal{P}_{\Gamma}^{\text{NI}}(\mathbf{a},\mathbf{b})$ illustrated in the previous section):
\[
(\mathcal{P}_{\Gamma}^{\text{NI}}(\mathbf{a},\mathbf{b}),\eta_1)
\overset{\Phi_\text{LGV}}{\sijto}(\mathcal{P}_{\Gamma}(\mathbf{a},\mathbf{b}),\eta_1)
\overset{\Psi}{\sijto} (\mathcal{P}_{\Gamma}(\mathbf{a},\mathbf{b}),\eta_2)
\overset{\Phi_\text{LGV}^{-1}}{\sijto} (\mathcal{P}_{\Gamma}^{\text{NI}}(\mathbf{a},\mathbf{b}),\eta_2).
\]
In fact, this sijection $\Psi$ is constructed by rotating each path 180 degrees around the midpoint of the two endpoints of the path.
It is important that we no longer need to handle the non-intersecting condition, thanks to the LGV lemma.

\section{Double LGV lemma technique} \label{sec::DLGV}
In this section, we formalize the method used in the bijective proof provided in the previous section.
Let $\mathbf{a}=(a_1,a_2,\ldots,a_n)$ and $\mathbf{b}=(b_1,b_2,\ldots,b_n)$ be vertices on a directed acyclic graph $\Gamma_1$, and $\mathbf{c}$ and $\mathbf{d}$ be vertices on a directed acyclic graph $\Gamma_2$.
According to the bijective proof of the LGV lemma, we have sijections $\mathcal{P}^{\text{NI}}_{\Gamma_1}(\mathbf{a},\mathbf{b}) \sijto \mathcal{P}_{\Gamma_1}(\mathbf{a},\mathbf{b})$ and
$\mathcal{P}^{\text{NI}}_{\Gamma_2}(\mathbf{c},\mathbf{d}) \sijto \mathcal{P}_{\Gamma_2}(\mathbf{c},\mathbf{d})$. Furthermore, these sijections are compatible with many statistics.
Therefore, by constructing a (good) sijection
$\mathcal{P}_{\Gamma_1}(\mathbf{a},\mathbf{b}) \sijto \mathcal{P}_{\Gamma_2}(\mathbf{c},\mathbf{d})$,
we can establish a property between combinatorial objects (e.g., plane partitions) associated with non-intersecting paths,
without explicitly handling the non-intersecting conditions.
This method is not entirely new; for example, it has been used in the bijective proof of the hook-content formula for the number of plane partitions \cite{RW}.

As a further application of this method, we provide another reason why a Schur polynomial $s_\lambda$, defined combinatorially as
\[
	s_\lambda(x_1,x_2,\ldots,x_n) = \sum_{T \colon \text{semi-standard Young tableaux of shape $\lambda$ with entries at most $n$}} x^T,
\]
where $\lambda$ is a partition, is symmetric.
Note that a bijective proof of this statement is already known and can be found in \cite[Theorem 7.10.2]{Sta2}.
Let $\Gamma$ be the directed acyclic graph defined in Section \ref{sec::main}, and we set $a_i = (-i,-i)$ and $b_j = (-j+\mu_j, -j+n-\mu_j)$, where $\mu=(\mu_1,\mu_2,\ldots,\mu_{\ell(\lambda)}) := \lambda^\text{t}$ is the transpose of $\lambda$. Then, as is well known, semi-standard Young tableaux of shape $\lambda$ with entries at most $n$ are in bijection to non-intersecting paths
$\mathcal{P}_\Gamma^{\text{NI}}(\mathbf{a},\mathbf{b})$, where each path records the data in the corresponding column.
Especially, an entry whose value is $t$ corresponds to an edge in $E_t := \{ (x,y) \to (x,y-1) \mid x-y = t-1 \} \cup \{ (x,y) \to (x+1,y) \mid x-y = t-1 \}$.
For a permutation $\sigma$ of $n$, let $\eta_\sigma = ( \#(P \cap E_{\sigma_1}), \#(P \cap E_{\sigma_2}), \ldots, \#(P \cap E_{\sigma_n}))$, where $\#(P \cap E_t)$ means the (multiplicity) number of edges in $E_t$ that constitute the paths $P$.
Now, proving the symmetry of the Schur function is equivalent to constructing a compatible sijection $(\mathcal{P}_\Gamma^{\text{NI}}(\mathbf{a},\mathbf{b}), \eta_\id) \sijto (\mathcal{P}_\Gamma^{\text{NI}}(\mathbf{a},\mathbf{b}), \eta_\sigma)$ for all permutation $\sigma$.
According to the method, i.e., by applying the LGV lemma two times, it is sufficient to construct a compatible sijection $(\mathcal{P}_\Gamma(\mathbf{a},\mathbf{b}), \eta_\id) \sijto (\mathcal{P}_\Gamma(\mathbf{a},\mathbf{b}), \eta_\sigma)$, and this can be carried out by permuting the directions (south or east) in each path according to $\sigma$.
Note that this sijection is actually a sign-preserving bijection.

\section{Conclusion}
In this paper, we provide a refined enumeration of plane partitions of a given shape with bounded entries (Theorem \ref{thm::main}), with a bijective proof.
We abstracted the method used in the bijective proof and, as an application, illustrate yet another reason why the Schur function is symmetric.
This method is an application of the notion of compatibility defined in \cite{Inoue}, and we believe both the notion of compatibility and this method can be applied to many bijective problems.

\section*{Acknowledgements}
I would like to thank my supervisor, Prof. Ralph Willox for helpful discussions and comments on the manuscript.
The author is partially supported by FoPM, WINGS Program, the University of Tokyo and JSPS KAKENHI No.\ 23KJ0795.

\end{document}